\documentclass[12pt]{article}
\usepackage{amssymb,amsmath,mathrsfs}
\usepackage{hyperref}
\usepackage{graphicx}
\usepackage{color}
\usepackage[OT2,OT1]{fontenc}
\usepackage{upquote}
\usepackage{authblk}
\usepackage{mathtools}

\usepackage{array,multirow}
\usepackage[numbers,sort&compress]{natbib}

\usepackage{float}

\begin{document}

\title{\Large\bf A two-domain MATLAB implementation \\
for efficient computation of the \\
Voigt/complex error function}

\date{October 5, 2022}

\author{
{\small Sanjar M. Abrarov, Rehan Siddiqui, Rajinder K. Jagpal} \\
\vspace{-0.4cm}
{\small and Brendan M. Quine}
}

\maketitle

\begin{abstract}
In this work we develop a new algorithm for the efficient computation of the Voigt/complex error function. In particular, in this approach we propose a two-domain scheme where the number of the interpolation grid-points is dependent on the input parameter $y$. The error analysis we performed shows that the MATLAB implementation meets the requirements for radiative transfer applications involving the HITRAN molecular spectroscopic database. The run-time test shows that this MATLAB implementation provides rapid computation, especially at smaller ranges of the parameter $x$.
\vspace{0.25cm}
\vspace{0.25cm}
\\
\end{abstract}

\section{Introduction}

The complex error function, also known as the Faddeeva function, is defined as~\cite{Faddeeva1961, Armstrong1967, Gautschi1970, Abramowitz1972}
\begin{equation}\label{eq_1}
w\left(z\right) = e^{-z^2}\left(1 + \frac{2i}{\sqrt\pi}\int\limits_0^z{e^{t^2}dt}\right),
\end{equation}
where $z = x + iy$ is a complex argument. The~complex error function $w\left( z \right)$ is closely related to the complex probability function~\cite{Armstrong1967}
\[
W\left(z\right) = \frac{i}{\pi}\int\limits_{-\infty}^\infty{\frac{e^{-t^2}}{z - t}}dt.
\]

The complex probability function can be written in terms of its real and imaginary parts~\cite{Armstrong1967}
$$
W\left(z\right) = K\left(x,y\right) + iL\left(x,y\right)
$$
such that
\begin{equation}\label{eq_2}
K\left(x,y\right) = \operatorname{Re}\left[W\left(z\right)\right] = \frac{y}{\pi}\int\limits_{-\infty}^\infty{\frac{e^{-t^2}}{y^2 + {\left(x - t\right)^2}}}dt
\end{equation}
and
\begin{equation}\label{eq_3}
L\left(x,y\right) = \operatorname{Im}\left[W\left(z\right)\right] = \frac{1}{\pi}\int\limits_{-\infty }^\infty{\frac{e^{-t^2}\left(x - t\right)}{y^2 + \left(x - t\right)^2}}dt.
\end{equation}

Both functions $w\left(z\right)$ and $W\left(z\right)$ are equal to each other on the upper half of the complex plane, when $y = \operatorname{Im}\left[z\right] > 0$ \cite{Armstrong1967}. Consequently, it follows that
$$
w\left(z\right) = K\left(x,y\right) + iL\left(x,y\right), \qquad y > 0.
$$
Further we will imply that the parameter $y = \operatorname{Im}\left[z\right]$ is always greater than~zero.

The real part of the complex error function $K\left(x,y\right)$ is known as the Voigt function~\cite{Faddeeva1961, Armstrong1967, Gautschi1970, Abramowitz1972} that is widely used in Atmospheric Science to describe emission and absorption of the photons by atmospheric molecules~\cite{Edwards1988, Quine2002, Jagpal2010, Berk2017, Pliutau2017, Siddiqui2017, Siddiqui2021, Jagpal2022}. Specifically, the~Voigt function is used to compute wavelength-dependent absorption coefficients by using the HITRAN molecular spectroscopic database~\cite{Hill2016}. The~imaginary part of the complex error function $L\left(x,y\right)$ is also used in many applications. For~example, it can describe the spectral behavior of the index of refraction in various materials~\cite{Balazs1969,Chan1986}.

Despite simple representations, the~integrals \eqref{eq_1}, \eqref{eq_2} and \eqref{eq_3} do not have analytical solutions and must be computed numerically. Many approximations are available in the scientific literature~\cite{Humlicek1982,Drummond1985,Poppe1990,Schreier1992,Kuntz1997,Wells1999,Letchworth2007,Abrarov2009,Pagnini2010,Imai2010,Zaghloul2011,Abrarov2011,Amamou2013,Abrarov2015}. Although~this problem has been known for many decades, the derivation of new approximations for the complex error function $w\left(z\right)$ and the development of their efficient algorithms still remain an interesting topic~\cite{Abrarov2018a,Abrarov2018b,Wang2019,Abrarov2019,Kumar2020,Nordebo2021,Pliutau2021,Azah2021}.

In our previous publication we proposed an algorithm for the efficient computation of the complex error function based on a single-domain implementation with vectorized interpolation~\cite{Abrarov2019}. However, despite rapid performance it has some limitations, including a limited range for the parameter $x$ as well as restrictions for the input array type. As~a further development, in~this work we present a new MATLAB implementation without those drawbacks. The~numerical analysis and computational tests we performed show that the proposed algorithmic implementation meets all the requirements in terms of accuracy and run-time performance for the efficient computation of the Voigt/complex error function in the radiative transfer~applications.

\section{Approximations}

\subsection{Sampling Based~Approximation}

In our {work}~\cite{Abrarov2015} we have proposed a new sampling method based on incomplete cosine expansion of the sinc function. In~particular, it is shown that, using a new sampling method based on incomplete cosine expansion of the sinc function, we can obtain the following approximation
\begin{equation}\label{eq_4}
w\left(z\right) \approx \sum\limits_{m = 1}^M {\frac{a_m + b_m\left(z + i\varsigma/2\right)}{c_m^2 -\left(z + i\varsigma/2\right)^2}},
\end{equation}
where the expansion coefficients are given by
$$
a_m = \frac{\sqrt\pi\left(m - 1/2\right)}{2M^2 h}\sum\limits_{n = - N}^N{e^{{\varsigma ^2}/4 - {n^2}{h^2}}\sin\left(\frac{\pi\left(m - 1/2\right)\left(nh + \varsigma/2\right)}{Mh}\right)},
$$
$$
b_m = -\frac{i}{M\sqrt{\pi}}\sum\limits_{n = -N}^N{e^{\varsigma^2/4 - n^2 h^2}\cos\left(\frac{\pi\left(m - 1/2\right)\left(nh + \varsigma/2\right)}{Mh}\right)},
$$
$$
c_m = \frac{\pi\left(m - 1/2\right)}{2M\pi},
$$
with parameters $N$, $M$, $h$ and $\varsigma$ that can be taken as $23$, $23$, $0.25$ and $2.75$, respectively.

For more rapid performance in algorithmic implementation, it is reasonable to define the function
$$
\Omega\left(z\right) = \sum\limits_{m = 1}^M{\frac{a_m + b_m z}{c_m^2 - z^2}},
$$
such that Equation \eqref{eq_4} can be represented as~\cite{Abrarov2018a}
\begin{equation}\label{eq_5}
w\left(z\right) \approx \Omega\left(z + i\varsigma/2\right).
\end{equation}

Overall, approximation \eqref{eq_5} is highly accurate. However, its accuracy deteriorates with decreasing parameter $y$. In~order to resolve this problem we can use the following approximation~\cite{Abrarov2018a}
\begin{equation}\label{eq_6}
w\left(z\right) \approx e^{-z^2} + z\sum\limits_{m = 1}^{M + 2}{\frac{\alpha_m - \beta_m z^2}{\gamma_m - \theta_m z^2 + z^4}},
\end{equation}
where the expansion coefficients are
$$
\alpha_m = b_m\left[c_m^2 - \left(\frac{\varsigma^2}{2}\right)^2\right] + i\alpha _m\varsigma,
$$
$$
\beta_m = b_m,
$$
$$
\gamma_m = c_m^4 + \frac{c_m^2\varsigma^2}{2} + \frac{\varsigma^4}{16}
$$
and
$$
\theta_m = 2c_m^2 - \frac{\varsigma^2}{2}.
$$
It is interesting to note that this approximation of the complex error function is obtained by substituting approximation \eqref{eq_4} into the right side of the identity (see~\cite{Abrarov2018a} for details in~derivation)
$$
w\left(z\right) = e^{-z^2} + \frac{w\left(z\right) - w\left(-z\right)}{2}.
$$

For $\left| z \right| > 8$ one of the best choices is the approximation based on the Laplace continued fraction~\cite{Faddeeva1961,Gautschi1970}. In~particular, in~our algorithm we used
\begin{equation}\label{eq_7}
w\left(z\right) \approx \frac{\left(i/\sqrt\pi\right)}{z - \frac{1/2}{z - \frac{1}{z - \frac{3/2}{z - \frac{2}{z - \frac{5/2}{z - \frac{3}{\frac{7/2}{z - \frac{4}{\frac{9/2}{z - \frac{5}{z - \frac{11/2}{z}}}}}}}}}}}}.
\end{equation}

This algorithm is implemented in MATLAB as a script file {\it fadsamp.m} that utilizes three approximations \eqref{eq_5}, \eqref{eq_6} and \eqref{eq_7} as follows~\cite{Abrarov2018a}
\begin{equation}\label{eq_8}
w\left(z\right) \approx \left\{
\begin{aligned}
{\rm  {Equation~\eqref{eq_5}}
}, &\quad{\rm if}\:(\left|x + iy\right| \le 8)\,\cap\,(y > 0.05x), \\
{\rm  {Equation~\eqref{eq_6}}}
, &\quad{\rm if}\:(\left|x + iy\right| \le 8)\,\cap\,(y \le 0.05x), \\
{\rm  {Equation~\eqref{eq_7}}}
, &\quad{\rm otherwise.}
\end{aligned}
\right.
\end{equation}
As we have shown in our publication~\cite{Abrarov2018a}, approximation \eqref{eq_8} provides highly accurate and rapid computation of the complex error function without poles and can be used to cover the entire complex~plane.

\subsection{Modified Trapezoidal~Rule}

In 1945, English mathematician and cryptanalyst Alan Turing, who succeeded to decrypt sophisticated  machine codes of the Enigma during the Second World War~\cite{Copeland2004}, published an interesting paper where he proposed an elegant method of numerical integration for some class of integrals~\cite{Turing1945}. Nowadays, his method of the numerical integration that involves the residue calculus is regarded as the modified trapezoidal rule~\cite{Trefethen2014} or the generalized trapezoidal rule~\cite{Azah2021}. The~comprehensive and detailed description of Turing's method of integration may be found in literature~\cite{Trefethen2014}.

In 1949, Goodwin showed how to implement Turing's idea to the integrals of kind~\cite{Goodwin1949}
$$
\int\limits_{-\infty}^\infty{f\left(x\right)e^{-x^2}dx}.
$$

Using the method described by Goodwin, Chiarella and Reichel in their work~\cite{Chiarella1968} derived the series expansion for the following integral
$$
\Psi\left(x,t\right) = U\left(x,t\right) + iV\left(x,t\right) = \frac{\Omega\left(x,t\right)}{\left(4\pi t\right)^{1/2}}\int\limits_{-\infty}^\infty{\frac{e^{-u^2}}{u^2 + \Omega^2\left(x,t\right)}du},
$$
where $\Omega\left(x,t\right) = \left(1 - ix\right)/\left(2t^{1/2}\right)$. In~particular, they showed that that the function $\Psi \left( {x,t} \right)$ can be approximated as a series
\begin{equation}\label{eq_9}
\begin{gathered}
\Psi\left(x,t\right) \approx \frac{h}{\Omega\left(x,t\right)\left(4\pi t\right)^{1/2}} + \frac{2h\Omega\left(x,t\right)}{\left(4\pi t\right)^{1/2}}\sum\limits_{n = 1}^\infty{\frac{e^{-n^2 h^2}}{\Omega ^2\left(x,t\right) + n^2 h^2}} \\
+ \frac{\pi e^{\Omega^2\left(x,t\right)}}{\left(\pi t\right)^{1/2}\left(1 - e^{2\pi\Omega\left(x,t\right)/h}\right)}H\left(t - \frac{h^2}{\pi^2}\right), \\
\end{gathered}
\end{equation}
where $h$ is a small fitting parameter and $H\left(t\right)$ is the Heaviside step function defined as
\[
H\left(t\right) = \left\{
\begin{aligned}
0, &\quad {\rm if} \: t < 0, \\
1/2, &\quad {\rm if} \: t = 0, \\
1, &\quad {\rm if} \: t > 0.
\end{aligned}
\right.
\]

Thus, due to Heaviside step function, \mbox{Equation \eqref{eq_9}} can be separated into three~parts
\begin{equation}\label{eq_10}
\Psi\left(x,t\right) \approx \frac{h}{\Omega\left(x,t\right)\left(4\pi t\right)^{1/2}} + \frac{2h\Omega\left(x,t\right)}{\left(4\pi t\right)^{1/2}}\sum\limits_{n = 1}^\infty{\frac{e^{-n^2 h^2}}{\Omega^2\left(x,t\right) + n^2 h^2}} , \quad t < \frac{h^2}{\pi^2},
\end{equation}
\begin{equation}\label{eq_11}
\begin{aligned}
\Psi\left(x,t\right) \approx &\frac{h}{\Omega\left(x,t\right)\left(4\pi t\right)^{1/2}} + \frac{2h\Omega\left(x,t\right)}{\left(4\pi t\right)^{1/2}}\sum\limits_{n = 1}^\infty{\frac{e^{-n^2 h^2}}{\Omega^2\left(x,t\right) + n^2 h^2}} \\
&+\frac{\pi e^{\Omega^2\left(x,t\right)}}{2\left(\pi t\right)^{1/2}\left(1 - e^{2\pi\Omega\left(x,t\right)/h}\right)}, \qquad t = \frac{h^2}{\pi^2},
\end{aligned}
\end{equation}
\begin{equation}\label{eq_12}
\begin{aligned}
\Psi\left(x,t\right) \approx &\frac{h}{\Omega\left(x,t\right)\left(4\pi t\right)^{1/2}} + \frac{2h\Omega\left(x,t\right)}{\left(4\pi t\right)^{1/2}}\sum\limits_{n = 1}^\infty{\frac{e^{-n^2 h^2}}{\Omega^2\left(x,t\right) + n^2 h^2}} \\
&+\frac{\pi e^{\Omega^2\left(x,t\right)}}{\left(\pi t\right)^{1/2}\left(1 - e^{2\pi\Omega\left(x,t\right)/h}\right)}, \qquad t > \frac{h^2}{\pi^2}.
\end{aligned}
\end{equation}

Equation \eqref{eq_11} deals only with a single point $h^2/\pi^2$ for the parameter $t$ and does not represent any practical interest. Therefore, further we will consider only two equations, Equations \eqref{eq_10} and \eqref{eq_12}.

Mata and Reichel~\cite{Matta1971} showed the relations
$$
K\left(x,y\right) = \frac{1}{y\sqrt\pi}U\left(\frac{x}{y},\frac{1}{4y^2}\right)
$$
and
$$
L\left(x,y\right) = \frac{1}{y\sqrt\pi}V\left(\frac{x}{y},\frac{1}{4y^2}\right)
$$
that link both functions $\Psi\left(x,t\right)$ and $w\left(x,y\right)$ with each other. Consequently, using these relations the series expansions \eqref{eq_10} and \eqref{eq_12} can be reformulated as
\begin{equation}\label{eq_13}
w\left(z\right) \approx \frac{2ihz}{\pi}\sum\limits_{k = 0}^N{\frac{e^{-t_k^2}}{z^2 - t_k^2}}
\end{equation}
and
\begin{equation}\label{eq_14}
w\left(z\right) \approx \frac{2e^{-z^2}}{1 + e^{-2i\pi z/h}} + \frac{2ihz}{\pi}\sum\limits_{k = 0}^N{\frac{e^{-t_k^2}}{z^2 - t_k^2}},
\end{equation}
respectively, where $t_k = \left(k + 1/2\right)h$ and $h$ can be chosen to be equal to $\sqrt{\pi/\left(N + 1 \right)}$ \cite{Azah2021}.

Consider an expansion series for the complementary error function that was reported by Hunter and Regan~\cite{Hunter1972} (see also~\cite{Azah2021})
$$
{\rm erfc}\left(z\right) \approx \frac{2hz e^{-z^2}}{\pi}\sum\limits_{k = 1}^N{\frac{e^{-\left(k - 1/2\right)^2h^2}}{z^2 + \left(k - 1/2\right)^2h^2}}  + \frac{2}{1 + e^{2\pi z/h}}, \qquad x < \frac{\pi}{h}.
$$
Using the identity relating complex error function and complementary error function~\cite{Gautschi1970}
$$
w\left(z\right) = e^{-z^2}{\rm erfc}\left(-iz\right),
$$
we obtain
\begin{equation}\label{eq_15}
w\left(z\right) \approx \frac{2e^{-z^2}}{1 + e^{-2i\pi z/h}} + \frac{ih}{\pi z} + \frac{2ihz}{\pi}\sum\limits_{k = 1}^N{\frac{e^{-\tau_k^2}}{z^2 - \tau_k^2}},
\end{equation}
where $\tau_k = kh$ \cite{Azah2021}.

Equations \eqref{eq_13} and \eqref{eq_14} have poles at $z = t_k$ while Equation \eqref{eq_15} has poles at $z = \tau_k$. Furthermore, each of these equations can cover with high accuracy only in the corresponding domain. However, Al Azah and Chandler-Wilde~\cite{Azah2021} recently proposed the following~approximation
\begin{equation}\label{eq_16}
w\left(z\right) \approx \left\{
\begin{aligned}
{\rm Equation~(13)}, &\qquad{\rm if}\: y \ge \max\left(\pi,x\right),\\
{\rm Equation~(15)}, &\qquad{\rm if}\: y < x \cap 1/4 \le\varphi\left(x/h\right)\le 3/4,\\
{\rm Equation~(14)}, &\qquad{\rm otherwise,}
\end{aligned}
\right.
\end{equation}
where $\varphi\left(t\right) = t - \left\lfloor t \right\rfloor\in\left[0,1\right)$, that appeared to be very efficient since it can be used for rapid and highly accurate computation without poles at $N = 11$. This is possible to achieve since, according to approximation \eqref{eq_16}, Equations \eqref{eq_13}, \eqref{eq_14} and \eqref{eq_15} are used interchangeably depending on the domain over the entire complex~plain.

\section{Algorithmic~Implementation}

Previously we have reported a new algorithm based on a vectorized interpolation over a single-domain~\cite{Abrarov2019}. Such an approach provides accuracy better than $10^{-6}$ at $y \ge 10^{-8}$ for the HITRAN~\cite{Hill2016} applications. However, this implementation has several limitations. In~particular, there is a limitation $\left|x\right| \le 10^5$. Although~it is possible to increase the range for more than $10^5$, it requires to introduce more {interpolation grid-points} for precomputation. Furthermore, this MATLAB implementation accepts an input only as a vector (one dimensional array) or as a~scalar.

One of the ways to implement an efficient algorithm is to use an approximation based on the two-domain scheme that we proposed in our earlier publication~\cite{Abrarov2009}
\[
\tag{17a}\label{eq_17a}
K\left( {x,y} \right) \approx \left\{
\begin{aligned}
&{\rm interpolation,} &&\qquad \frac{x^2}{27^2} + \frac{y^2}{15^2} \le~1 \\
&\frac{a_1 + b_1 x^2}{a_2 + b_2 x^2 + x^4}, &&\qquad \frac{x^2}{27^2} + \frac{y^2}{15^2} > 1,
\end{aligned}\right.
\]
where the coefficients are~\cite{Kuntz1997}
\[
\begin{aligned}
&a_1 = y/\left(2\sqrt{\pi}\right)+y^3/\sqrt{\pi}\approx 0.2820948y + 0.5641896y^3 \\
&b_1 = y/\sqrt{\pi}\approx~0.5641896y \\
&a_2 = 0.25 + y^2 + y^4 \\
&b_2 =-1 + 2y^2
\end{aligned}
\]
such that
$$
\frac{a_1 + b_1 x^2}{a_2 + b_2 x^2 + x^4}=\operatorname{Re}\left\{\frac{\left(i/\sqrt{\pi}\right)}{z-\frac{1/2}{z}}\right\}.
$$
Consequently, the~complex error function can also be approximated as~\cite{Abrarov2019}
\[
\tag{17b}\label{eq_17b}
w\left(z\right) = K\left(x,y\right)+iL\left(x,y\right) \approx \left\{
\begin{aligned}
&{\rm interpolation,} &&\qquad \frac{x^2}{27^2} + \frac{y^2}{15^2} \leqslant~1 \\
&\frac{\left(i/\sqrt{\pi}\right)}{z-\frac{1/2}{z}}, &&\qquad \frac{x^2}{27^2} + \frac{y^2}{15^2} > 1.
\end{aligned}
\right.
\]

Although the algorithms for the Voigt and complex error functions, built on Equations \eqref{eq_17a} and \eqref{eq_17b} can provide rapid computations, their accuracies deteriorate at~$y < 10^{-6}$.

In order to resolve this problem we developed a new algorithm that utilizes the internal MATLAB built-up features. Unlike interpolation algorithms shown in~\cite{Abrarov2009,Abrarov2019}, the~proposed approach implies that the number of the interpolation grid-points required for precomputation is not constant and dependent on the input parameter $y$ such that
\begin{equation}\label{eq_18}\tag{18}
N_{gp} = \frac{1}{\sqrt{y}} + \delta,
\end{equation}
where $\delta$ is the offset that was found experimentally to be $5 \times 10^3$. The~corresponding interpolation grid-points range is bounded by $\left[-r,r\right]$, where $r = 35$ is the radius that was also found experimentally. The~interpolation grid-points are distributed non-equidistantly in logarithmic scale to increase their density near the origin along the $x$ axis.

Equation \eqref{eq_18} is entirely empirical. It is reasonable to suggest that the number of the interpolation grid-points should be increased with decreasing $y$. However, taking the reciprocal of $y^n$, such that $n \ge 1$, results in a very rapid increase of the interpolation grid-points $N_{gp}$ with decreasing $y$. Our experimental results show that a more or less optimal value $N_{gp}$ can be found by taking the reciprocal of $\sqrt{y}$. The~parameter $\delta$ in {Equation} \eqref{eq_18} is to provide a sufficient number of the interpolation grid-points $N_{gp}$ within the internal domain $|x+iy| \le r$ as parameter $y$ increases.

The use of Equation \eqref{eq_18} is given as a basic for applications. If~higher accuracy is required, the~number of the interpolation grid-points $N_{gp}$ can be further increased
by using, for~example, a~modified form of the equation above
$$
N_{gp} = \frac{2}{\sqrt{y}} + 3\delta.
$$
Such a modification improves accuracy by an order of the magnitude.

The algorithm utilizes two domains, internal and external, that are bounded by a circle of radius $\left|x + iy\right| = r$. The internal domain is situated within the circle while the external domain is situated outside~of it.

All points within internal domain $\left|x + iy\right| \le r$ are computed by the MATLAB built-in interpolation function {\it interp1} though {interpolation grid-points} that are computed by using the function file {\it fadsamp.m} provided in our article~\cite{Abrarov2018a}. The spline method was found to be the best for~interpolation.

All points outside external domain $\left|x+iy\right| > r$ are computed by the following approximation
$$
w\left(z\right) \approx \frac{\left(i/\sqrt \pi\right)}{z - \frac{1/2}{z - \frac{1}{z - \frac{3/2}{z}}}}
$$
that represents a simplified version of Equation \eqref{eq_7} above.

The parameter $y$ is dependent on temperature and pressure. Therefore, it is related to the atmospheric height. In~a radiative transfer model, the atmospheric column is sliced into layers \cite{Edwards1988, Quine2002}. The layers can be considered homogeneous if their thicknesses are small enough. Typically, a thickness of $50\:m$ (or $20$ layers per $km$) is sufficient to consider the layer homogeneous for atmospheric modeling. This implies that, within a given layer, the~temperature and pressure are constants. Since the density of air in the mesosphere is extremely low, its contribution to the radiance is practically negligible. Consequently, taking into consideration that only the troposphere and stratosphere, with an atmospheric column corresponding to a height of up to 40--50\:km, actually contribute to the radiance, we can conclude that the total number of layers as well as elements for parameter $y$ does not exceed $1000$ per gas in atmosphere. However, in~many practical tasks, for a single value $y$ at a given atmospheric temperature and pressure, we may need to cover a wide spectral range with high resolution. This signifies that the number of elements for parameter $x$ can exceed a million. Since the number of elements for parameter $y$ is much smaller than the number of elements for parameter $x$, it is reasonable to implement a configuration in which parameter $x$ is an array and parameter $y$ is a scalar. This approach is more rapid in the MATLAB implementation since the application of a 2D array $M \times N$, where $N$ and $M$ are numbers of elements for parameters $x$ and $y$, respectively, requires considerably more computer~memory.

As it has been mentioned above, the~computation of the {interpolation grid-points} in its original version of the function file {\it w2dom.m} is performed by external function {\it fadsamp.m}. However, any other MATLAB function file that can provide highly accurate computation of the complex error function may also be used for computation {of interpolation grid-points}. For~example, the~function files such as {\it fadf.m} \cite{Abrarov2011} and {\it fadfunc.m} \cite{Abrarov2018b} can also be used as alternatives. The~script of the function file {\it w2dom.m} is given in Appendix A.

\section{Error~Analysis}

In order to exclude the rounding and truncation errors in computation by using the most recent HITRAN database~\cite{Hill2016}, the~values of the $K\left(x,y\right)$ and $L\left(x,y\right)$ functions with $6$ or more accurate decimal digits in their mantissas are required. Therefore, in~radiative transfer applications involving the HITRAN database, the~accuracy of computation of these functions has to be better than $10^{-6}$.

\begin{figure}[H]
\begin{center}
\includegraphics[width=28pc]{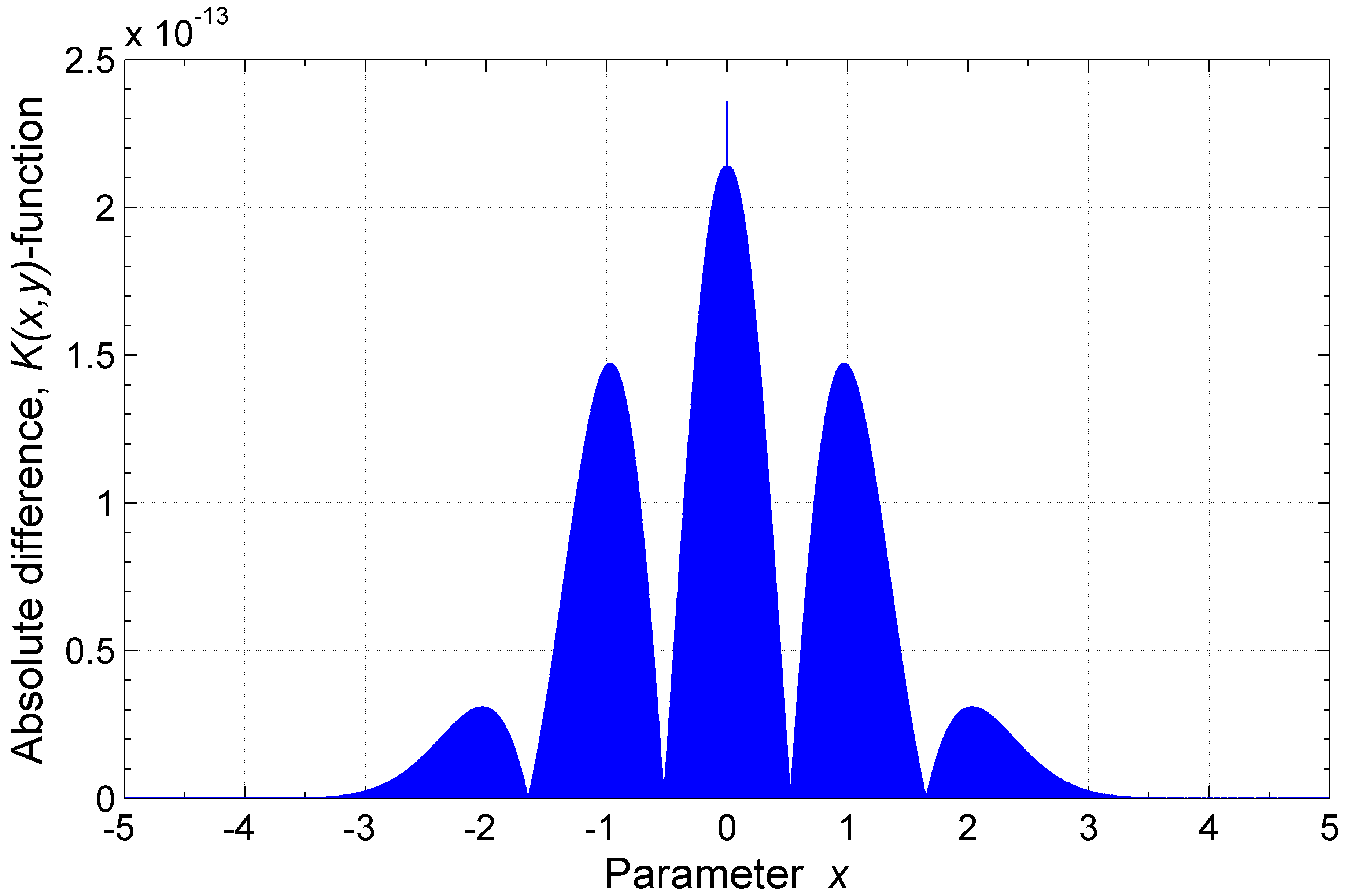}\hspace{2pc}%
\begin{minipage}[b]{24pc}
\vspace{0.3cm}
{\sffamily {\bf{Fig. 1.}} Absolute difference for the real part $K\left(x,y\right)$ of the complex error function.}
\end{minipage}
\end{center}
\end{figure}

Figure~1 shows the absolute difference $\left|K_{ref.}\left(x,y\right) - K\left(x,y\right)\right|$, where $K_{ref.}\left(x,y\right)$ is a highly accurate reference, in~the range $-5 \le x \le 5$ at $y = 10^{-8}$. As~we can see, the~absolute difference does not exceed $2.5 \times 10^{-13}$.

Figure~2 shows the absolute difference $\left|L_{ref.}\left(x,y\right) - L\left(x,y\right)\right|$, where $L_{ref.}\left(x,y\right)$ is a highly accurate reference, in~the range $-5 \le x \le 5$ also at $y = 10^{-8}$. We can see that for the imaginary part the absolute difference also does not exceed $2.5 \times 10^{-13}$.

\begin{figure}[H]
\begin{center}
\includegraphics[width=28pc]{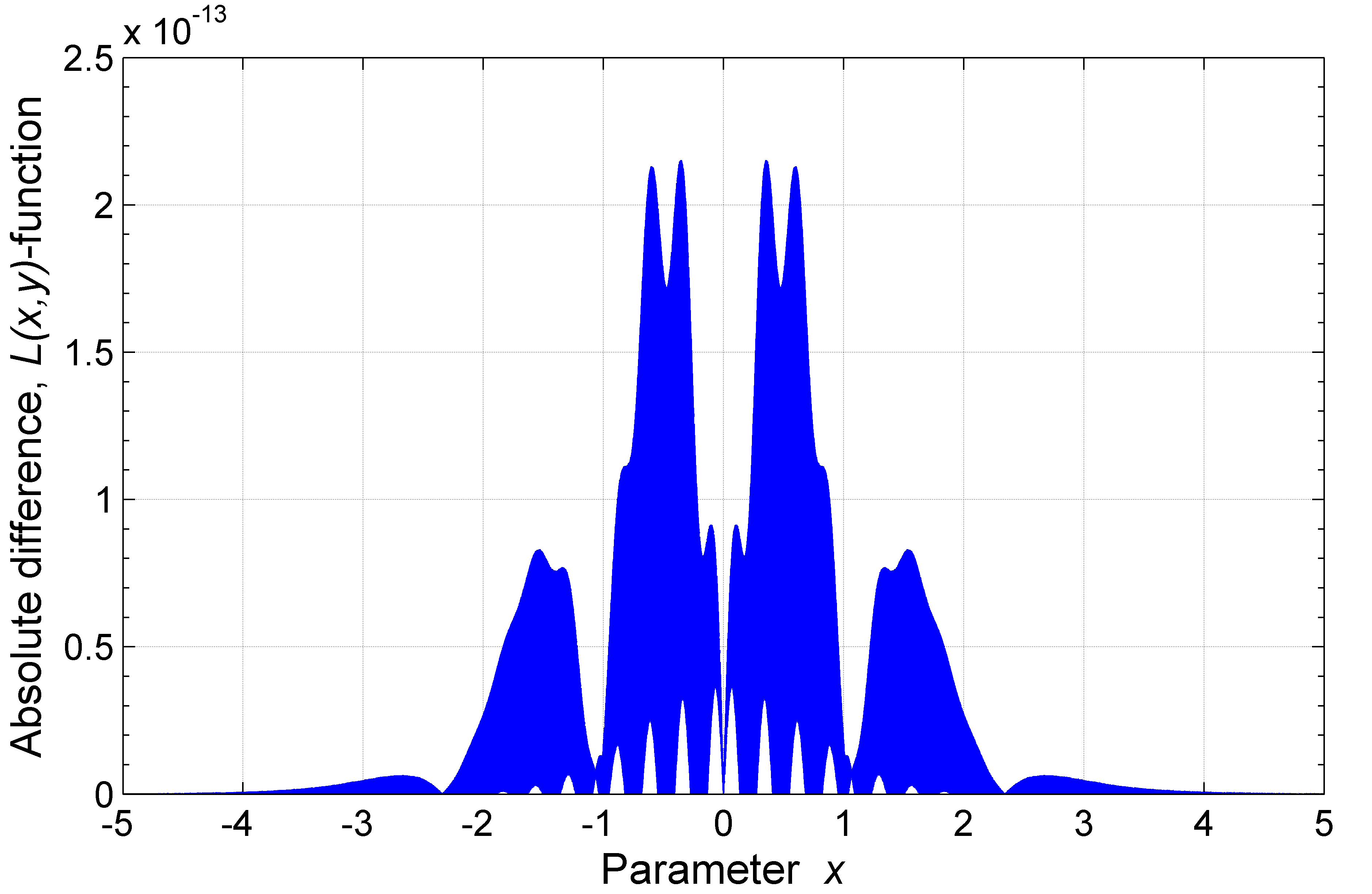}\hspace{2pc}%
\begin{minipage}[b]{24pc}
\vspace{0.3cm}
{\sffamily {\bf{Fig. 2.}} Absolute difference for the imaginary part $L\left(x,y\right)$ of the complex error function.}
\end{minipage}
\end{center}
\end{figure}

For more rigorous error analysis, we can apply the following relative errors
$$
\Delta_{\operatorname{Re}} = \frac{\left|K_{ref.}\left(x,y\right) - K\left(x,y\right)\right|}{K_{ref.}\left(x,y\right)}
$$
and
$$
\Delta_{\operatorname{Im}} = \frac{\left|L_{ref.}\left(x,y\right) - L\left(x,y\right)\right|}{L_{ref.}\left(x,y\right)}
$$
for the real and imaginary parts of the complex error function $w\left(z\right)$, respectively.

Figure 3 depicts the relative error for the real part of the complex error function $\operatorname{Re} \left[z\right] = K\left(x,y\right)$ in the domain $0 \le x \le 50$ and $0 \le y \le 10^{1.699}$, where  $10^{1.699} \approx 50$. Inset in Figure~3 shows the smaller domain near the origin $0 \le x \le 5.5$ and $0 \le y \le 10^{0.74}$, where $10^{0.74} \approx 5.5$. Beyond~this smaller domain the observed color is blue indicating that the relative error is below $10^{-12}$. As~we can see from Figure~3, the~relative error does not exceed $\sim$$10^{-10}$.

\begin{figure}[H]
\begin{center}
\includegraphics[width=28pc]{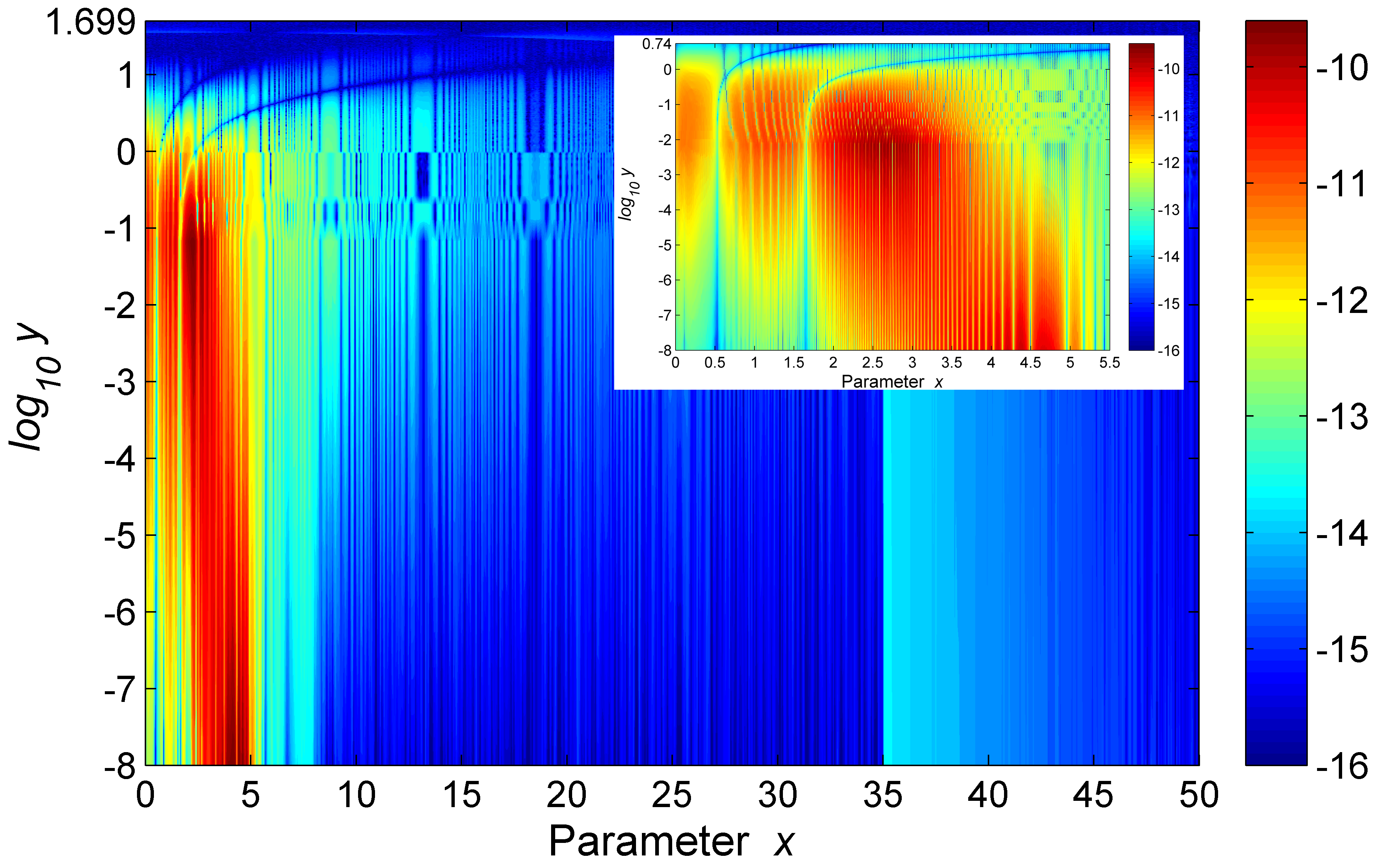}\hspace{2pc}%
\begin{minipage}[b]{24pc}
\vspace{0.3cm}
{\sffamily {\bf{Fig. 3.}} Relative error for the real part $K\left(x,y\right)$ of the complex error function.}
\end{minipage}
\end{center}
\end{figure}

\begin{figure}[H]
\begin{center}
\includegraphics[width=28pc]{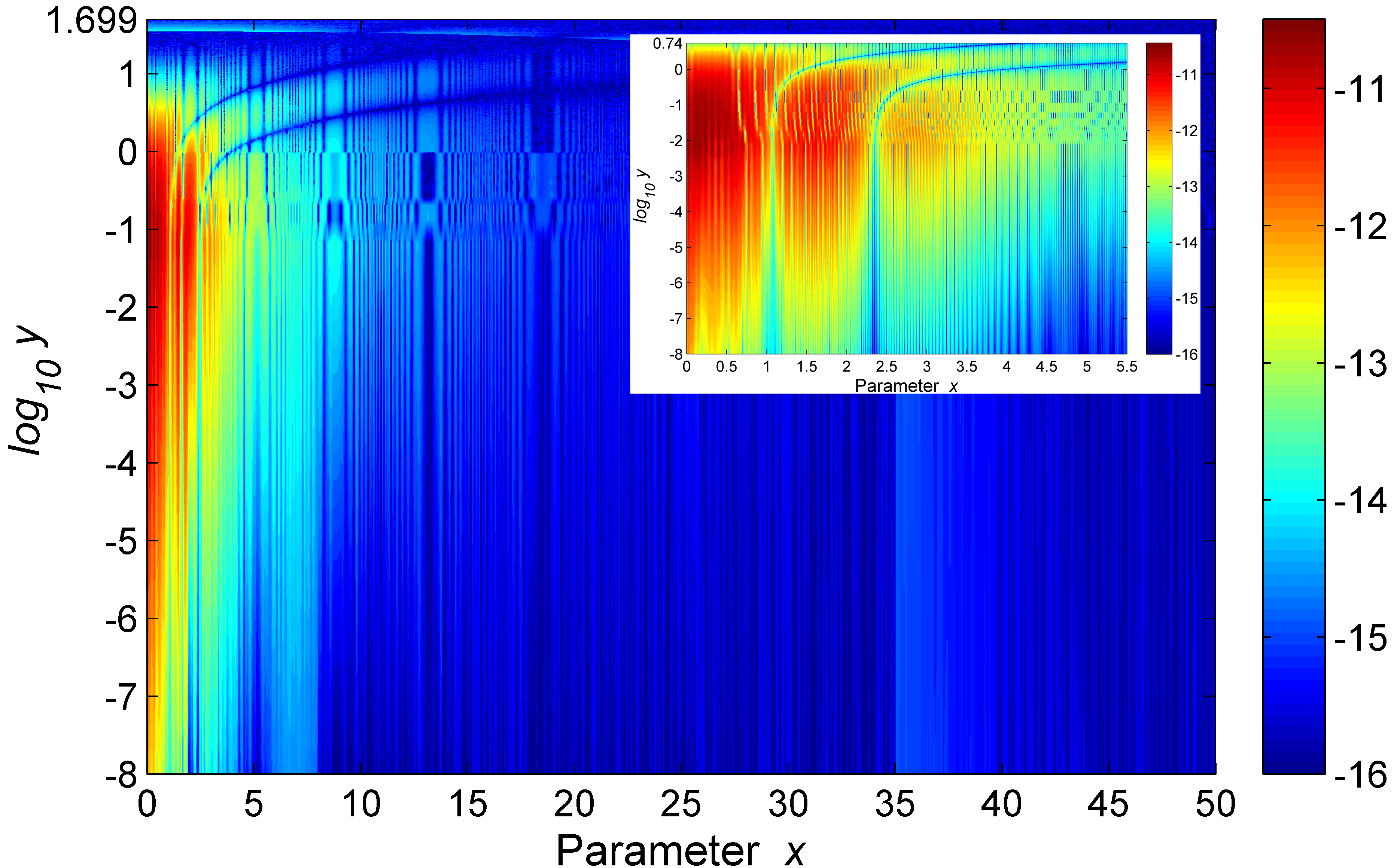}\hspace{2pc}%
\begin{minipage}[b]{24pc}
\vspace{0.3cm}
{\sffamily {\bf{Fig. 4.}} Relative error for the imaginary part $L\left(x,y\right)$ of the complex error function.}
\end{minipage}
\end{center}
\end{figure}

Figure 4 illustrates the relative error for the imaginary part of the complex error function $\operatorname{Im}\left[z\right] = L\left(x,y\right)$ in the domain $0 \le x \le 50$ and $0 \le y \le 10^{1.699}$. Inset in Figure~4 shows the smaller domain  near origin $0 \le x \le 5.5$ and $0 \le y \le 10^{0.74}$. Beyond~this smaller domain the color is blue and the relative error is below $10^{-12}$. As~we can see from Figure~4, the~relative error is generally lower in the imaginary part and does not exceed $\sim 10^{-11}$.

It is commonly known that the accuracy of $K\left(x,y\right)$ and $L\left(x,y\right)$ functions tend to deteriorate when $y$ decreases~\cite{Armstrong1967,Humlicek1982,Wells1999,Letchworth2007,Amamou2013}. However, from Figures~3 and 4 we can see that the decrease of the parameter $y$ does not deteriorate the accuracy. This is possible to achieve since in accordance with Equation \eqref{eq_18} the number of the {interpolation grid-points} $N_{gp}$ increases with decreasing $y$. Therefore, this technique enables us to resolve efficiently this problem in computation of the complex error function $w\left(z\right)$.

\section{Run-Time~Test}

The run-time test was performed with MATLAB function files {\it wTrap.m} \cite{Azah2021}, {\it fadsamp.m}~\cite{Abrarov2018b} and {\it w2dom.m} {at equidistantly distributed $10$ million points for parameter $x$ at $y = 10^{-8}$.} The results of the run-time test are shown in Table~1.

\begin{table}[H]
\begin{center}
\begin{tabular}{|c|c|c|c|}
\hline
\multirow{2}{*}{\bf Algorithm} & \multicolumn{3}{c|}{\bf Run-time in seconds}
\\
\cline{2-4}
& $x\in\left[-10,10\right]$
& $x\in\left[-10^2,10^2\right]$
& $x\in\left[-10^3,10^3\right]$
\\
\cline{1-4}
{\it wTrap.m}   & 2.41 & 2.55 & 2.45 \\
{\it fadsamp.m} & 4.14 & 1.78 & 1.54 \\
{\it w2dom.m}   & 1.23 & 1.04 & 0.98 \\
\hline
\end{tabular}
\caption{Run-time of algorithms for $10$ million points at three different ranges.}
\end{center}
\end{table}

It has been reported that both algorithms {\it wTrap.m} and {\it fadsamp.m} are highly accurate in computation~\cite{Azah2021}. In~particular, the~maximum values in relative errors for the algorithms {\it wTrap.m} and {\it fadsamp.m} are found to be $\sim$$10^{-15}$ and $\sim$$10^{-14}$, respectively. However, the~computational speed of these two algorithms differs depending on the range for the input parameter $x$. In~particular, Table~1 shows that at smaller range for parameter $x$ the function file  {\it wTrap.m} performs computation faster than the function file {\it fadsamp.m}. However, as~the range for  parameter $x$ increases, our algorithm {\it fadsamp.m} becomes faster. The~run-time test also reveals that the algorithm {\it w2dom.m} is always faster regardless of the chosen range. This is particularly evident when the range for the input parameter $x$ is~smaller.

Table~2 shows the results of the run-time test at $30$ million equidistantly distributed points {for parameter $x$ at $y = 10^{-8}$.} As we can see, the~algorithm {\it w2dom.m} remains more rapid proportionally at extended size of the input array as compared to the {\it wTrap.m} and {\it fadsamp.m} algorithms.

\begin{table}[H]
\begin{center}
\begin{tabular}{|c|c|c|c|}
\hline
\multirow{2}{*}{\bf Algorithm} & \multicolumn{3}{c|}{\bf Run-time in seconds}
\\
\cline{2-4}
& $x\in\left[-10,10\right]$
& $x\in\left[-10^2,10^2\right]$
& $x\in\left[-10^3,10^3\right]$
\\
\cline{1-4}
{\it wTrap}   & 8.38 & 8.33 & 8.42 \\
{\it fadsamp} & 14.37 & 6.02 & 5.21 \\
{\it w2dom}   & 3.78  & 3.07 & 2.86 \\
\hline
\end{tabular}
\caption{Run-time of algorithms for $30$ million points at three different ranges.}
\end{center}
\end{table}

The MATLAB code for the run-time test is provided in Appendix B. The~scripts of function files {\it wTrap.m} and {\it fadsamp.m} can be accessed from the cited literature~\cite{Azah2021,Abrarov2018a}, respectively.

\section{Conclusions}

A new algorithm for the efficient computation of the Voigt/complex error function is shown. We propose a two-domain scheme where the number of {interpolation grid-points} $N_{gp}$ is dependent on parameter $y$. The~error analysis shows that our MATLAB implementation meets the requirements for radiative transfer applications utilizing the HITRAN {molecular} spectroscopic database. The~run-time test we performed shows that this MATLAB implementation can provide rapid computation especially at smaller ranges of parameter~$x$.

\section*{Acknowledgments}

This work is supported by the National Research Council Canada, Thoth Technology Inc., York University, Epic College of Technology and Epic Climate Green (ECG) Inc.

\pagestyle{empty}
\section*{Appendix A}
\vspace{-0.5cm}
\footnotesize
\begin{verbatim}

function FF = w2dom(x,y,opt)

% SYNOPSIS:
%          x   - array or scalar
%          y   - scalar
%          opt - option for the real and imaginary parts
%
% opt = 1 returns value(s) for the Voigt function, K(x,y)-function
% opt = 2 returns value(s) for the L(x,y)-function
% opt = 3 returns value(s) for the complex error function
%
% This code is primarily intended to work in the range y >= 1e-8 for the
% HITRAN applications for accelerated computation of the Voigt/complex
% error function.
%
% -------------------------------------------------------------------------
% Example:
%          x = linspace(-10,10,1e7); y = 1e-8; tic; w2dom(x,y); toc
% -------------------------------------------------------------------------
if ~isscalar(y)
    disp('Input parameter y must be a scalar')
    return % terminate computation if y is not a scalar
end

if nargin == 2
    opt = 3;
    disp('Default value opt = 3 is assigned.')
end

if opt ~= 1 && opt ~= 2 && opt ~= 3
    disp(['Wrong parameter opt = ',num2str(opt),'! Use 1, 2 or 3.'])
    return
end

if y < 1e-8
    if opt == 1
        FF = real(fadsamp(x+1i*y));
    elseif opt == 2
        FF = imag(fadsamp(x+1i*y));
    else
        FF = fadsamp(x+1i*y);
    end
    return
end

FF = zeros(size(x));
radius = 35; % define radius
ind = abs(x + 1i*y) <= radius;

gp =[-radius,radius]; % define grid-points
if ~isempty(x(ind))

    offset = 5*1e3; % assign offset
    nump = 1/sqrt(y) + offset; % assign number of points
    % For better accuracy use for example nump = 2/sqrt(y) + 3*offset;
    
    gp = radius*(logspace(log10(1 + eps),log10(2),nump) ...
    - 1); % notice the log scale

    gp = [flip(-gp),gp];
end

switch opt
    case 1
        FF(ind) = real(internD(x(ind),y,gp));
        FF(~ind) = real(externD(x(~ind) + 1i*y));
    case 2
        FF(ind) = imag(internD(x(ind),y,gp));
        FF(~ind) = imag(externD(x(~ind) + 1i*y));
    otherwise
        FF(ind) = internD(x(ind),y,gp);
        FF(~ind) = externD(x(~ind) + 1i*y);
end

    function IntD = internD(x,y,gp) % internal domain
        IntD = interp1(gp,fadsamp(gp + 1i*y),x, ...
            'spline'); % interpolated values
    end

    function ExtD = externD(z) % external domain

        num = 1:4; % define a row vector
        num = num/2;

        ExtD = num(end)./z; % start computing from the end
        for m = 1:length(num) - 1
            ExtD = num(end - m)./(z - ExtD);
        end
        ExtD = 1i/sqrt(pi)./(z - ExtD);
    end
end

\end{verbatim}
\normalsize

\section*{Appendix B}
\vspace{-0.5cm}
\footnotesize
\begin{verbatim}

clear
clc

% Table for run-time in seconds
tab{1,1}='Algorithm';
tab{1,2}='-10 to 10'; % range 1
tab{1,3}='-10^2 to 10^2'; % range 2
tab{1,4}='-10^3 to 10^3'; % range~3

tab{2,1}='wTrap';
tab{3,1}='fadsamp';
tab{4,1}='w2dom';

y=1e-8; % smallest value for the parameter y
maxN=10; % max number of~cycles

for k=[1,3] % k is a factor
    for m=1:3 % three ranges
        x0=10^m; x=linspace(-x0,x0,k*1e7); % 1 and 3 are for 1e7 and 3*1e7 ...
                                           % points, respectively

        if k==1
            disp(['10 million points in range ',num2str(m)])
        else
            disp(['30 million points in range ',num2str(m)])
        end
        disp('Computing, please wait!')

        tic; for n=1:maxN; wTrap(x+1i*y,11); end; tab{2,m+1}=toc/maxN;
        tic; for n=1:maxN; fadsamp(x+1i*y); end; tab{3,m+1}=toc/maxN;     
        tic; for n=1:maxN; w2dom(x,y); end; tab{4,m+1}=toc/maxN;

        clc
    end
    if k==1; tab1=tab; else tab2=tab; end % assign tab1 or tab2
end

disp('Displaing Table 1:')
disp(tab1)

disp('Displaing Table 2:')
disp(tab2)

\end{verbatim}
\normalsize

\pagestyle{plain}

\end{document}